\documentclass[a4paper,11pt]{article}
\usepackage{inputenc}
\usepackage{graphicx}
\usepackage{amssymb}
\usepackage{latexsym, bm}
\usepackage{multicol}
\usepackage{indentfirst}
\usepackage{amssymb,amsfonts}
\usepackage{amsmath}
\usepackage{setspace}
\usepackage{enumerate}
\usepackage{CJK}
\newtheorem{theorem}{Theorem}

\newtheorem{lemma}{Lemma}
\newtheorem{corollary}{Corollary}

\newtheorem{conjecture}{Conjecture}

\usepackage{amssymb}
\begin{spacing}{1.0}

\begin{document}
\date{}

\title{Ramsey numbers of large books\footnote{Center for Discrete Mathematics, Fuzhou University,
Fuzhou, 350108, P.~R.~China.  Email: {\em chenxjq@foxmail.com, linqizhong@fzu.edu.cn, chunlin\_you@163.com.}  }}

\author{Xun Chen,\;\;
Qizhong Lin,\footnote{Corresponding author. Supported in part by NSFC(No.\ 1217010182) and ``New century excellent talents support plan for institutions of higher learning of Fujian province"(SJ2017-29).} \;\;
Chunlin You
}
\maketitle
\begin{abstract}
A book $B_n$ is a graph which consists  of $n$ triangles sharing a common edge.
In 1978, Rousseau and Sheehan conjectured that the Ramsey number satisfies $r(B_m,B_n)\le 2(m+n)+c$ for some constant $c>0$.
In this paper, we obtain that $r(B_m, B_n)\le 2(m+n)+o(n)$ for all $m\le n$ and $n$ large, which confirms the conjecture of Rousseau and Sheehan asymptotically. As a corollary, our result implies that a related conjecture of Faudree, Rousseau and Sheehan (1982) on strongly regular graph holds asymptotically.
\medskip

\end{abstract}

\section{Introduction}
For graphs $H_1$ and $H_2$, the Ramsey number $r(H_1,H_2)$ is defined as the smallest integer $N$ such that for any blue/red edge coloring of $K_N$, there exists either a blue $H_1$ or a red $H_2$.

For $k\ge2$, let $B_n^{(k)}$ be the book graph that consists of $n$ copies of $K_{k+1}$ sharing a $K_k$ in common. When $k=2$, we write $B_n$ instead of $B_n^{(2)}$ for convenience. The Ramsey number of books has been studied extensively.
Erd\H{o}s et al. \cite{efrs}, and later Thomason \cite{tho} by using the constructive method obtained that $r(B_n^{(k)},B_n^{(k)})\ge(2^k+o_k(1))n.$
Conlon \cite{Conlon} proved that $r(B_n^{(k)},B_n^{(k)})=(2^k+o_k(1))n$, which confirms a conjecture of Thomason \cite{tho} asymptotically and also gives an answer to a problem proposed by Erd\H{o}s \cite{efrs}. Using a different method, the upper bound has been improved to $r(B_n^{(k)},B_n^{(k)})\le2^kn+O_k(\frac{n}{(\log \log\log n)^{1/25}})$ by Conlon, Fox and Wigderson \cite{c-f}. For more Ramsey numbers of large books versus other graphs, we refer the reader to \cite{lr,lp,nr,N-R-S,sud} and other related references.

Another seminal result on book Ramsey numbers by Rousseau and Sheehan \cite{R-S} is as follows, in which the general upper bound was also pointed out by Parson \cite{Par}.
\begin{theorem}\label{r-s-g}
 We have $r(B_1,B_n)=2n+3$ for $n\ge2$.
Moreover, if $2(m+n)+1>(n-m)^2 /3$, then $$r(B_m, B_n)\le 2(m+n+1).$$
Generally, $r(B_m, B_n)\le m+n+2+\lfloor \frac{2}{3} \sqrt{3(m^2+mn+n^2)}\rfloor.$
\end{theorem}

The authors made the following conjecture,
one can also see \cite{F-R-S,N-R1}.
\begin{conjecture}\label{r-s}
There exists a constant $c>0$ such that for all $m,n\ge1$,
$$r(B_m,B_n)\le 2(m+n)+c.$$
\end{conjecture}

Theorem \ref{r-s-g} implies that the above conjecture holds when $m$ and $n$ are nearly equal.
Since Rousseau and Sheehan \cite{R-S} obtained the first exact value that $r(B_m,B_n)=2n+3$ when $m=1$ and $n\ge2$, researchers are interested in finding more exact values of $r(B_m,B_n)$. 
 Faudree, Rousseau and Sheehan \cite{F-R-S} proved that $r(B_m,B_n)=2n+3$ for all $m\ge2$ and $n\ge (m-1)(16m^3+16m^2-24m-10)+1$. In \cite{F-R-S}, the authors also proved that $r(B_2,B_n)=2n+6$ for $n=2,5,11$. In \cite{N-R1}, Nikiforov and Rousseau proved that $r(B_m,B_n)=2n+3$ for all $n\ge 10^6m$. Using a stability result on books by Bollob\'{a}s and Nikiforov \cite{b-n}, Nikiforov and Rousseau \cite{N-R2} further proved the following result.
 \begin{theorem}\label{nr}
If $m\le n/6-o(n)$ and $n$ is large, then $$r(B_m,B_n)=2n+3.$$
Moreover, the constant $1/6$ is asymptotically best possible.
\end{theorem}


Now let us turn to the {\em strongly regular graph} which was introduced by Bose \cite{bs}. A strongly regular graph $srg(\nu,k,\lambda,\mu)$ is a graph of order $\nu$ which is $k$-regular, in which any pair of vertices have $\lambda$ common neighbors if they are adjacent, and $\mu$ common neighbors otherwise. It is clear that its complement is also a strongly regular graph with parameters $\nu$, $\nu-k-1$, $\nu-2k+\mu-2$, and $\nu-2k+\lambda.$

Strongly regular graphs provide good lower bounds for the related Ramsey numbers.
It is well known that for each prime power $q=4n+1$, the Paley graph $P_q$ is a strongly regular graph $srg(q,\frac{q-1}{2},\frac{q-5}{4},\frac{q-1}{4})$ which is self-complementary.
Paley graphs give optimal lower bounds for $r(K_3,K_3)$ and $r(K_4,K_4)$. Almost all known lower bounds for $r(K_n,K_n)$ with $n$ small are obtained from Paley graphs, see \cite{mathon,shearer}.
In particular, Rousseau and Sheehan \cite{R-S} applied Paley graphs to give that $r(B_n,B_n)=4n+2$ when $4n+1$ is a prime power.
In \cite{b-s}, Bose and Shrikhande constructed strongly regular graphs $srg(4k^2-1,2k^2,k^2,k^2)$ for all $k$ of the form  $k= 3^m2^{m+n-1}$ where $m,n\ge0$, not both zero, which together with Theorem \ref{r-s-g} yields that $r(B_{k^2-2},B_{k^2+1})=4k^2$ with such restriction.  For more constructions of strongly regular graphs and exact values of $r(B_m,B_n)$, we refer the reader to \cite{bro,N-R1}. In \cite{F-R-S}, Faudree, Rousseau and Sheehan also proposed the following conjecture, which is closely related to Conjecture \ref{r-s}.
 \begin{conjecture}\label{f-r-s}
There exists a constant $c\ge2$ such that for every $srg(\nu,k,\lambda,\mu)$,
$$2(s+t)-\nu\le c,$$
where $s=k-\lambda-1$ and $t=k-\mu$.
\end{conjecture}
{\em Remark.} As pointed in \cite{F-R-S}, Conjecture \ref{r-s} indeed implies Conjecture \ref{f-r-s}.

\medskip

In this paper, we confirm Conjecture \ref{r-s} by Rousseau and Sheehan \cite{R-S} asymptotically.

\begin{theorem}\label{New1} For all $m\le n$ and $n$ large,
$$r(B_m, B_n)\le 2(m+n)+o(n).$$
\end{theorem}

 Let us point out that Theorem \ref{nr} due to Nikiforov and Rousseau \cite{N-R2} already implies that Conjecture \ref{r-s} holds for $m\le n/6-o(n)$ and large $n$.
Set $m=\lfloor\alpha n\rfloor$. From the general upper bound 
 in Theorem \ref{r-s-g}, we get
$$r(B_m, B_n)\le \left(1+\alpha+\frac23\sqrt{3(\alpha^2+\alpha+1)}\right)n+3.$$
The coefficient of the upper bound in Theorem \ref{New1} is $2+2\alpha+o(1)$, which is clearly less than $1+\alpha+\frac23\sqrt{3(\alpha^2+\alpha+1)}$ for any $0<\alpha<1$ and large $n$.
Therefore, our result improves the general upper bound 
 in Theorem \ref{r-s-g} for any $0<\alpha<1$ and large $n$.

\medskip

As a corollary, we confirm Conjecture \ref{f-r-s} asymptotically.
\begin{corollary}
For every $srg(\nu,k,\lambda,\mu)$ and large $\nu$,
$$2(s+t)-\nu=o(\nu),$$
where $s=k-\lambda-1$ and $t=k-\mu$.
\end{corollary}
{\bf Proof.} Suppose that the assertion fails for some strongly regular graph $G=srg(\nu,k,\lambda,\mu)$ with $\nu$ large, i.e., $2(s+t)-\nu>\epsilon_0\nu$  for some $\epsilon_0>0$. Note that the complement of $G$ is a strongly regular graph
$$srg(\nu,\nu-k-1,\nu-2k+\mu-2,\nu-2k+\lambda).$$
Let $m=\lambda+1$ and $n=\nu-2k+\mu-1$. Since $G$ contains no $B_m$ by noting any pair of vertices have $\lambda$ common neighbors if they are adjacent in $G$ and similarly its complement contains no $B_n$, and
\begin{align*}
\nu-2(m+n)&=\nu-2(\lambda+1)-2(\nu-2k+\mu-1)
\\&=2[(k-\lambda-1)+(k-\mu)]-\nu+2
\\&=2(s+t)-\nu+2
\\&>\epsilon_0\nu,
\end{align*}
it follows that
\[
r(B_m,B_n)>\nu>2(m+n)+\epsilon_0\nu,
\]
which contradicts Theorem \ref{New1}.
\hfill$\Box$

\medskip
In this paper, we also give a lower bound for $r(B_m, B_n)$ as follows.
\begin{theorem}\label{New2}
For any fixed $0< \alpha\le 1$, $r(B_{\lceil\alpha n\rceil}, B_n)\ge (\sqrt{4\alpha}+1+\alpha-o(1))n$.
\end{theorem}

 The above lower bound is asymptotically best possible when $\alpha=1$ from the above Theorem \ref{New1} or Theorem \ref{r-s-g} by Rousseau and Sheehan \cite{R-S}, and one can also see a much more general result established by Conlon \cite{Conlon}. For the special case when $\alpha=\frac{1}{2}$, the result above together with Theorem \ref{New1} yields that
  $$2.9142n< r(B_{\lceil n/2\rceil}, B_n)\le(3+o(1))n,$$ and the gaps between the lower and upper bounds from Theorems \ref{New1}--\ref{New2} will become smaller and smaller as $\alpha$ increases.


\section{Preliminaries}

Let $G=(V,E)$ be a graph. For two sets $A,B\subseteq V(G)$ (not necessarily disjoint), we write $e_G(A,B)$ for the number of edges between $A$ and $B$ in $G$, where each edge in $A\cap B$ will be counted twice. We call $$d_G(A,B)=\frac{e_G(A,B)}{|A||B|}$$ the density of the pair $(A,B)$. If $A=\{u\}$, then we will write $d_G(u,B)$ instead of $d_G(\{u\},B)$ for convenience. We always delete the subscription when there is no confusion.

Given $\epsilon>0$, a pair $(A,B)$ is called \emph{$\epsilon$-regular} if
$|d(A,B)-d(X,Y)|\le \epsilon$ for every $X\subseteq A$, $Y\subseteq B$ with $|X|\ge \epsilon |A|$ and $|Y|\ge \epsilon |B|$.
We say that a subset $U$ of the vertex set of a graph $G$ is $\epsilon$-regular if the pair $(U,U)$ is $\epsilon$-regular.
For graph $G$, we say a partition $V(G)=\cup_{i=1}^k V_i$ of $G$ is equitable if $||V_i|-|V_j||\le 1$ for all distinct $i$ and $j$.


A refined version of the regularity lemma \cite[Lemma 3]{Conlon} guarantees that one can find a regular subset in each part of the partition for any graph.
A key ingredient of the proof of our main result is the following refined regularity lemma \cite[Lemma 2.1]{c-f}, which is a slight strengthening of that due to Conlon \cite{Conlon} and the usual version of Szemer\'{e}di's regularity lemma \cite{sze78}. An earlier refined version of the regularity lemma due to Alon, Fischer, Krivelevich and Szegedy \cite{afk} has been used in the proof of the induced removal lemma. For many applications, we refer the reader to \cite{K-S} and other related references.

\begin{lemma}\label{regular}
For every $\epsilon>0$ and $M_0\in \mathbb{N}$, there is some $M=M(\epsilon, M_0)>M_0$ such that for every graph $G$, there is an equitable partition $V(G)=\cup_{i=1}^k V_i$ into $M_0\le k\le M$ parts so that the following hold:

1. Each part $V_i$ is $\epsilon$-regular.

2. For every $1\le i\le k$, there are at most $\epsilon k$ values $1\le j\le k$ such that the pair $(V_i, V_j)$ is not $\epsilon$-regular.
\end{lemma}


We will use the following version of the counting lemma, which is similar to Nikiforov and Rousseau \cite[Corollary 11]{N-R2} in which it requires all $U_i$ to be different, and one can see Conlon \cite[Lemma 5]{Conlon} for a more general version. For the general local counting lemma, see R\"{o}dl and Schacht \cite[Theorem 18]{rodl-s}. 

\begin{lemma}\label{counting}Let $l\ge1$ be an integer and $0<\epsilon\le1/(l+1)$.
If $U_1, U_2, U_{3},\ldots, U_{l+2}$ are (not necessarily distinct) vertex sets with
$(U_1,U_2)$ $\epsilon$-regular, and $(U_i,U_{j})$ $\epsilon$-regular of density $d_{i,j}$ for all $1\le i\le 2 <j\le l+2$, then there is an edge between $U_1$ and $U_2$ which is contained in at least $\sum_{j=3}^{l+2}(d_{1,j}d_{2,j}-2\epsilon)|U_{j}|$
triangles with the third vertex in $\cup_{j=3}^{l+2} U_{j}$.
\end{lemma}

\medskip

We will frequently use the following consequence, which can be used to count extensions of cliques and thus estimate the size of books, see \cite[Corollary 2.6]{c-f}.
\begin{lemma}\label{extend}
Let $\eta,\delta\in (0,1)$ be parameters with $\eta\le \delta^3/k^2$. Suppose $U_1,\ldots, U_k$ are (not necessarily distinct) vertex sets in a graph $G$ and all pairs $(U_i,U_{j})$ are $\eta$-regular with $\prod_{1\le i<j\le k}d(U_i,U_j)\ge \delta$. Let $Q$ be a randomly chosen copy of $K_k$ with one vertex in each $U_i$ with $1\le i\le k$ and say that a vertex $u$ extends $Q$ if $u$ is adjacent to every vertex of $Q$. Then, for any $u$,
$$\Pr(u\mbox{ extends } Q)\ge \prod\limits_{i=1}^k d(u,U_i)-4\delta.$$
\end{lemma}


We also need the following simple fact.

\begin{lemma}\label{ineq1}
If $0\le x_1, x_2\le 1$, then
$(1-x_1)^2+(1-x_2)^2+2x_1x_2\ge 1.$
\end{lemma}
{\bf Proof.} It suffices to show that
$
1-2x_1+x_1^2-2x_2+x_2^2+2x_1x_2\ge0,
$
which is clear since the left equals $(1-x_1-x_2)^2$.\hfill$\Box$

\section{Proof of Theorem \ref{New1}}

Let $m=\lfloor\alpha n\rfloor$ where $0<\alpha\le1$ is a constant.
For any sufficiently small $\gamma$ with $0<\gamma<1/10$, let $N=\lceil(2+2\alpha+\gamma)n\rceil$. Consider a red/blue edge coloring of the complete graph on vertex set $[N]$, where $[N]=\{1,2,\dots,N\}$.
Let $R$ and $B$ denote the graphs induced by all red and blue edges, respectively. Set
\begin{align}\label{pms}
\delta=\frac{\alpha \gamma}{100},\ \ \epsilon=\frac{\delta^3}{4},\ \ \text{and}\ \ M_0=1/\epsilon.
 \end{align}


We begin by applying Lemma \ref{regular} to the red graph, with the parameters $\epsilon$ and $M_0$ as above.
We obtain an equitable partition $V=V(K_N)=\cup_{i=1}^k V_i$ with a bounded number of parts such that, for each $i$, $V_i$ is $\epsilon$-regular and there are at most $\epsilon k$ values $1\le j\le k$ such that the pair $(V_i, V_j)$ is not $\epsilon$-regular.
Note that the colors are complementary, the same assertion holds for the blue graph.
In the following, we assume that  $|V_i|=N/k$ where $1\le i\le k$ for convenience.
For large $N$, we lose very little by assuming this.

 Now we construct a reduced graph $F$ on vertex set $\{v_1,\ldots,v_k\}$, in which $v_i$ is adjacent to $v_j$   if and only if $(V_i, V_j)$ is $\epsilon$-regular.
For each vertex $v_i$, we assign a color $c_i$ to $v_i$ such that $c_i$ is red if $d_R(V_i)\ge1/2$ where $d_R(V_i)$ is the density of the red subgraph induced by $V_i$. Otherwise $c_i$ is blue.
By the pigeonhole principle, at least $k'\ge k/2$ such $c_i$'s are the same, say color $\mathcal{A}$, where $\mathcal{A}\in\{red, blue\}$.

\medskip
{\bf Case 1.} \ $\mathcal{A}$ is red.

\medskip
By relabelling if necessary, we may assume that $v_1,v_2,\ldots,v_{k'}$ are these red vertices. Let $F'$ be the subgraph of $F$ induced by the red vertices $v_i$ with $1\le i\le k'$.
For this case, we color the edge $v_i v_j$ red in $F'$ if  $d_R(V_i,V_j)\ge1-\delta$ and blue otherwise.
Since $v_1$ has at most $\epsilon k\le 2\epsilon k'$ non-neighbors coming from irregular pairs, it follows that there are at least $(1-2\epsilon)k'$ parts $V_j$ with $1\le j\le k'$ such that $(V_1,V_j)$ is $\epsilon$-regular. Let $J$ be the set of all these indices $j$ and let $U=\bigcup_{j\in J}V_j$ be the union of all these $V_j$.

Suppose that $F'$ contains only red edges.
 Since $d_R(V_i)\ge1/2$, we apply Lemma \ref{counting} with $U_1=U_2=V_1$ and $U_{2+j}=V_{b_j}$ for $b_j\in J$ to obtain that there exists a red edge in $V_1$ which is contained in at least
\begin{align*}\sum\limits_{b_j\in J} (d_R(V_1,V_{b_j})^2-2\epsilon)|V_{b_j}|\ge ((1-\delta)^2-2\epsilon)|J|\frac{N}{k}\ge(1-2\delta)\left(\frac{1}{2}-\epsilon\right)(2+2\alpha+\gamma)n\ge n
\end{align*}
red triangles by noting (\ref{pms}). Thus we obtain a red book $B_n$ as desired.

In the following, we may assume that there exists a blue edge in this reduced graph $F'$. This corresponds two parts, say $V_1,V_2$, satisfying $d_R(V_1),d_R(V_2)\ge \frac{1}{2}$ but $d_B(V_1,V_2)\ge \delta$.

For any vertex $u\in V$ in the original graph $K_N$, let $x_i(u)=d_B(u,V_i)$, $i=1,2$. By Lemma \ref{ineq1}, we know that
$x_1(u)x_2(u)+\frac{1}{2}\sum_{i=1,2}(1-x_i(u))^2\ge \frac{1}{2}.$
Summing this inequality over all vertices of $V$, we get that
$$\sum_{u\in V}x_1(u)x_2(u)+\frac{1}{2}\sum_{i=1,2}\sum_{u\in V}(1-x_i(u))^2\ge \frac{1}{2}N.$$
It follows that either
\begin{align}\label{s-p1}
\sum_{u\in V}x_1(u)x_2(u)\ge \frac{\alpha}{2+2\alpha}N,
 \end{align}
 or
\begin{align}\label{s-p2}
\frac{1}{2}\sum_{i=1,2}\sum_{u\in V}(1-x_i(u))^2\ge \frac{1}{2+2\alpha}N.
\end{align}

First, suppose that $\sum_{u\in V}x_1(u)x_2(u)\ge \frac{\alpha}{2+2\alpha}N$.
Let $(u_1,u_2)$ be a randomly chosen blue edge with $u_1\in V_1$ and $u_2\in V_2$.
We apply Lemma \ref{extend} with $U_1=V_1$ and $U_2=V_2$ to obtain that for any vertex $u\in V$,
\[
\Pr(u \;\; \text{extends} \;\; (u_1,u_2))\ge x_1(u)x_2(u)-4\delta.
\]
Thus, we have that the expected number of blue extensions of a randomly chosen blue edge spanned by $(V_1,V_2)$ is at least
$$\sum_{u\in V}(x_1(u)x_2(u)-4\delta)\ge \left(\frac{\alpha}{2+2\alpha}-4\delta\right)N=\left(\frac{\alpha}{2+2\alpha}-4\delta\right)(2+2\alpha+\gamma)n\ge \alpha n,$$
where the last inequality holds by noting (\ref{pms}). Therefore, there must exist some blue edge with at least $\alpha n$ blue
extensions, giving a blue $B_m$ as desired.

On the other hand, suppose that
$\frac{1}{2}\sum_{i=1,2}\sum_{u\in V}(1-x_i(u))^2\ge \frac{1}{2+2\alpha}N$.
Without loss of generality, we may assume that $\sum_{u\in V}(1-x_1(u))^2\ge \frac{1}{2+2\alpha}N$.
By a similar argument by applying Lemma \ref{extend} with $U_1=U_2=V_1$, we obtain that the expected number of extensions of a random red edge in $V_1$ is at least
$$\sum_{u\in V}\left((1-x_1(u))^2-4\delta\right)\ge \left(\frac{1}{2+2\alpha}-4\delta\right)N\ge n$$
by a similar computation as above. Therefore, we see that a randomly chosen red edge inside $V_1$ will
have at least $n$ red extensions in expectation. Hence, there must exist a red $B_n$.

\medskip
{\bf Case 2.} \ $\mathcal{A}$ is blue.

\medskip
The proof for this case is similar to the one for Case 1. We swap the red and the blue colors, $d_R(U,V)$ and $d_B(U,V)$.
We mention the main steps of the proof for this case as follows.

For this case, there are at least $k'\ge k/2$ vertices, say $v_1,v_2,\ldots,v_{k'}$ , that are colored blue since $\mathcal{A}$ is blue.
 Let $F'$ be the subgraph of $F$ induced by the blue vertices $v_i$ with $1\le i\le k'$.
 We color the edge $v_i v_j$ blue if $d_B(V_i,V_j)\ge1-\delta$ and red otherwise.
 There are at least $(1-2\epsilon)k'$ parts $V_j$ with $1\le j\le k'$ such that $(V_1,V_j)$ is $\epsilon$-regular. Let $J$ be the set of all these indices $j$ and let $U=\bigcup_{j\in J}V_j$ be the union of all of these $V_j$.

 If this reduced graph $F'$ has only blue edges, then we can easily get a blue $B_m$ by using Lemma \ref{counting}. So we may suppose that there is a red edge in this reduced graph $F'$. This corresponds two parts, say $V_1,V_2$, satisfying $d_B(V_1),d_B(V_2)\ge \frac{1}{2}$ but $d_R(V_1,V_2)\ge \delta$.

For any vertex $u$ in the original graph $K_N$, let $x_i(u)=d_B(u,V_i)$, $i=1,2$.
Similarly, Lemma \ref{ineq1} implies that
 either $\sum_{u\in V}x_1(u)x_2(u)\ge \frac{1}{2+2\alpha}N$, or
$\frac{1}{2}\sum_{i=1,2}\sum_{u\in V}(1-x_i(u))^2\ge \frac{\alpha}{2+2\alpha}N$.
Again, a similar argument as Case 1, there exists a red $B_n$ in the first case or a blue $B_m$ in the second case. 

The proof of Theorem \ref{New1} is complete now.
 \hfill$\Box$

\section{Proof of Theorem \ref{New2}}

 For $0< \alpha\le 1$ and sufficiently small $0<\eta<1/10$,  let $\beta=\sqrt{4\alpha-\eta}+(1+\alpha)$ and $N=\beta n$. Clearly, $\beta\le 4$.
 Color the edges of $K_N$ by red and blue independently and randomly, where each edge is colored red with probability $p=\sqrt{{1}/{\beta}}-\delta$ and blue with probability $q=1-p$, where $\delta=\eta/100$. Clearly, $\delta<\sqrt{1/\beta}$.
Let $\mathcal{E}_1$ and $\mathcal{E}_2$ be the following events:
\medskip

$\mathcal{E}_1$: there exists a red $B_n$;

$\mathcal{E}_2$: there exists a blue $B_{\lceil\alpha n\rceil}$.
\medskip

We aim to show that $\Pr(\mathcal{E}_1\cup\mathcal{E}_2)<1$ for all large $n$.
If $uv$ is red, then the number of common red neighbors of $u$ and $v$, denoted by $X$, is a random variable that has a binomial distribution with parameters $N-2$ and $p^2$. It is clear that
$E(X)=(N-2)(\sqrt{1/\beta}-\delta)^2\le (1/\beta-\sqrt{1/\beta}\delta)N$
by noting $\delta<\sqrt{1/\beta}$. Therefore, by Chernoff bound (see e.g. \cite{A-S}), we obtain that
\begin{align*}
\Pr(\mathcal{E}_1)\le {N\choose 2}p\times\Pr\left(X\ge N/\beta\right)
\le {N\choose 2}\Pr\left(X\ge E(X)+\sqrt{\frac{1}{\beta}}\delta N\right)< N^2 e^{-\delta^2N/2}=o(1).
\end{align*}

Now we bound $\Pr(\mathcal{E}_2)$ as follows.
If $uv$ is blue, then the number of common blue neighbors of $u$ and $v$, denoted by $Y$, is a random variable that has a binomial distribution with parameters $N-2$ and $q^2$. Note that
$E(Y)=(N-2)q^2\le N(1-\sqrt{1/\beta}+\delta)^2\le ((1-\sqrt{1/\beta})^2+2\delta)N$
again by noting $\delta<\sqrt{1/\beta}$.

\medskip
{\bf Claim.} $E(Y)\le (\frac{\alpha}{\beta}-\frac{1}{2}\delta)N$.

\medskip
{\bf Proof.} Let us denote
$\lambda=\frac{\alpha}{\beta}-((1-\sqrt{1/\beta})^2+2\delta).$
It suffices to show $\lambda\ge \frac{1}{2}\delta$.
Indeed,
\begin{align}\label{delta}
\lambda&=\frac{\alpha}{\beta}-\left(1-2\sqrt{\frac{1}{\beta}}+\frac{1}{\beta}+2\delta\right)
=\frac{1}{\beta}[2\sqrt{\beta}-(\beta+1-\alpha+2\beta\delta)]\nonumber
\\&=\frac{4\beta-(\beta+1-\alpha+2\beta\delta)^2}{\beta[2\sqrt{\beta}+(\beta+1-\alpha+2\beta\delta)]}.
\end{align}
Note that $\beta=\sqrt{4\alpha-\eta}+(1+\alpha)$, so $\beta\le4$ and $\beta+1-\alpha=\sqrt{4\alpha-\eta}+2\le4$.
Note also from the assumption that $\delta= \eta/100$ and $\eta$ is small, we obtain that $\beta[2\sqrt{\beta}+(\beta+1-\alpha+2\beta\delta)]\le 40$, and
\begin{align*}
4\beta-(\beta+1-\alpha+2\beta\delta)^2&=4\beta-(\beta+1-\alpha)^2-4\beta(\beta+1-\alpha)\delta-4\beta^2\delta^2
\\&=\eta-4\beta(\beta+1-\alpha)\delta-4\beta^2\delta^2
\\ &\ge 20\delta.
\end{align*}
Therefore, $\lambda\ge \frac{1}{2}\delta$. The claim follows as desired. \hfill$\Box$

\medskip

Now, from the above claim and applying the Chernoff bound again, we obtain that
\begin{align*}
\Pr(\mathcal{E}_2)\le {N\choose 2}q\times\Pr\left(Y\ge \frac{\alpha}{\beta}N\right)=o(1).
\end{align*}
Consequently, we have $\Pr(\mathcal{E}_1\cup\mathcal{E}_2)=o(1)$, which completes the proof.\hfill$\Box$

\bigskip
\noindent
{\bf Acknowledgement:} The authors would like to thank the anonymous referees for their valuable comments which greatly improve the presentation of this paper, in particular, the presentation of the proof for Theorem \ref{New1}.

\end{spacing}
\end{document}